\documentclass[12pt,a4paper]{article}
\usepackage{amssymb, amscd, amsthm, amsmath,latexsym, amstext,mathrsfs}
\usepackage[all]{xy}

\newcommand{\lra}{\longrightarrow}

\newcommand{\veps}{\varepsilon}

\newcommand{\z}{\zeta}

\newcommand{\wt}[1]{\widetilde{#1}}

\newcommand{\matr}[4]{\begin{pmatrix} {#1}& {#2}\\ {#3}&{#4}\end{pmatrix}}

\newcommand{\C}{\mathbb C}

\newcommand{\Z}{\mathbb Z}

\def\mod{\mathrm{mod}}

\DeclareMathOperator{\GL}{\mathrm GL}

\numberwithin{equation}{section}

\newtheorem{lemma}{Lemma}[section]
\newtheorem{cor}[lemma]{Corollary}
\newtheorem{prop}[lemma]{Proposition}
\newtheorem{thm}[lemma]{Theorem}

\theoremstyle{definition}

\theoremstyle{remark}

\begin{document}
\title{Irreducible representations and Artin $L$-functions of
quasi-cyclotomic fields}
\author{Sunghan Bae$^{(1)}$, Yong Hu$^{(2)}$, Linsheng Yin$^{(2)}$\\
(1) Department of Mathematics, KAIST, Daejeon, Korea\\
(2) Dept. of Math. Sciences, Tsinghua Univ., Beijing 100084, China}
\maketitle
\begin{abstract} We determine all irreducible representations of
primary quasi-cyclotomic fields in this paper. The methods can be applied
to determine the irreducible representations of any quasi-cyclotomic field.
We also compute the Artin $L$-functions for a class of quasi-cyclotomic fields.
\end{abstract}

\section{Introduction}
A quadratic extension of a cyclotomic field, which is non-abelian Galois
over the rational number field $\Bbb Q$, is called a quasi-cyclotomic field.
All quasi-cyclotomic fields are described explicitly in [8] followed the works
in [1] and [3]. They are generated by a canonical $\Z/2\Z$-basis. The minimal
quasi-cyclotomic field containing the quadratic roots of one element of the
basis is called a primary quasi-cyclotomic field.
L.Yin and C.Zhang [7] have studied the arithmetic of any quasi-cyclotomic field. In
this paper we determine all irreducible representations of primary quasi-cyclotomic
fields. The methods apply to determine the irreducible representations
of an arbitrary quasi-cyclotomic field. We also compute
the Artin $L$-functions for a class of quasi-cyclotomic fields.

First we recall the constructions of primary quasi-cyclotomic fields.
Let $S$ be the set consisting of $-1$ and all prime numbers. For $p\in S$,
we put $\bar p=4,8,p$ and set $p^*=-1, 2,(-1)^{\frac{p-1}2}p$ if $p=-1,2$
and an odd prime number, respectively. Let
$K=\Bbb Q(\zeta_{\bar pq})$ be the cyclotomic field of conductor $\bar pq$.
For a class $[a]\in\Bbb Q/\Bbb Z$, we set $\sin[a]=2\sin a\pi$ for $0<a<1$ and
$\sin[0]=1$. For prime numbers $p<q$, we define
\[
v_{pq} =
\prod_{i=0}^{\frac{p-1}2}\prod_{j=0}^{\frac{q-1}2}\frac{\sin[\frac{iq+j}{pq}]}
 {\sin[\frac{jp+i}{pq}]}\quad \qquad(p>2)
 \]
and
\[
v_{2q} =\frac{\sin[\frac 14]}{\sin[\frac 1{4q}]}\prod_{j=0}^{\frac{q-1}2}
\frac{\sin[\frac{j}{2q}]\cdot\sin[\frac{2j-1}{4q}]}
 {\sin[\frac{4j+i}{4q}]\cdot\sin[\frac{j}{q}]\cdot\sin[\frac{2j-1}{2q}]}.
\]
For $p<q\in S$, we put
\begin{displaymath}
u_{pq}:= \left \{
\begin{array}{ll}
\sqrt{q^*} &\text{if }p=-1\\
v_{pq} &\text{if }p=2\text{ or } p \equiv q \equiv 1 \mod 4 \\
\sqrt{p}\cdot v_{pq} &\text{if } p \equiv 1, ~ q \equiv 3 \mod 4 \\
\sqrt{q}\cdot v_{pq} &\text{if } p \equiv 3, ~ q \equiv 1 \mod 4 \\
\sqrt{pq}\cdot v_{pq} &\text{if } p \equiv q \equiv 3 \mod 4.
\end{array}
\right.
\end{displaymath}

Let $\tilde K=K(\sqrt{u_{pq}})$. Then $\tilde K$ is the minimal one in all quasi-cyclotomic
fields which contain $\sqrt{u_{pq}}$. We call these fields $\tilde K$ primary quasi-cyclotomic
fields. Let $G=\text{Gal}(K/\Bbb Q)$ and $\tilde G=\text{Gal}(\tilde K/\Bbb Q)$. We always
denote by $\epsilon$ the unique non-trivial element of $\text{Gal}(\tilde K/K)$.
If $(p,q)=(-1,2)$, then the group $G$ is generated by two elements $\sigma_{-1}$ and
$\sigma_2$, where $\sigma_{-1}(\zeta_8)=\zeta_8^{-1}$ and $\sigma_2(\zeta_8)=\zeta_8^{5}$.
If $p>2$, then $G$ is generated by two elements $\sigma_p$ and $\sigma_q$,
where $\sigma_p(\zeta_{p})=\zeta_{p}^a,\sigma_p(\zeta_{q})=\zeta_{q}$
and $\sigma_q(\zeta_{p})=\zeta_{p},\sigma_q(\zeta_{q})=\zeta_{q}^b$,
with $a,b$ being generators of $(\Bbb Z/p\Bbb Z)^*$ and $(\Bbb Z/q\Bbb Z)^*$
respectively. If $p=2$, then $G$ is generated by three elements $\sigma_{-1},\sigma_2$
and $\sigma_q$, where $\sigma_{-1},\sigma_2$ act on $\zeta_8$ as above and on $\zeta_q$
trivially, and $\sigma_q$ acts on $\zeta_q$ as above and on $\zeta_8$ trivially.

Next we describe the group $\tilde G$ by generators and relations.
An element $\sigma\in G$ has two liftings in $\tilde G$. By [Sect.3, 7] the action of the
two liftings on $\sqrt{u_{pq}}$ has the form $\pm\alpha\sqrt{u_{pq}}$ or
$\pm\alpha\sqrt{u_{pq}}/\sqrt{-1}$ with $\alpha>0$. We fix the lifting $\tilde\sigma$
of $\sigma$ to be the one with the positive sign. Then the other lifting of $\sigma$ is
$\tilde\sigma\epsilon$. The group $\tilde G$ is generated by $\epsilon,\tilde\sigma_p$
and $\tilde\sigma_q$ (and $\tilde\sigma_{-1}$ if $p=2$). Clearly $\epsilon$ commutes
with the other generators. In addition, we have
$\tilde\sigma_p\tilde\sigma_q=\tilde\sigma_q\tilde\sigma_p\epsilon$
(and $\tilde\sigma_{-1}$ commutes with $\tilde\sigma_2$ and $\tilde\sigma_q$ if $p=2$).
For an element $g$ of a group, we denote by $|g|$ the order of $g$ in the group.
Let $\log_{-1}:\{\pm 1\}\rightarrow\Bbb Z/2\Bbb Z$ be the unique isomorphism.
For an odd prime number $p$ and an integer $a$ with $p\nmid a$, let $(\frac ap)$ be the quadratic
residue symbol. We also define $(\frac a2)=(\frac a{-1})=1$ for any $a$. Then we have, see [7, Th.3],
\[
|\tilde\sigma_p|=(1+\log_{-1}(\frac{q^*}p))|\sigma_p|\quad
\text{ and }\quad |\tilde\sigma_q|=(1+\log_{-1}(\frac{p^*}q))|\sigma_q|,
\]
with the exception that $\tilde\sigma_2=2|\sigma_2|$ when $(p,q)=(-1,2)$. If $p=2$, we
have furthermore $|\tilde\sigma_{-1}|=|\sigma_{-1}|$.
Thus we have determined the group $\tilde G$ by generators and relations.

\section{Abelian subgroup of index 2}
In this section we construct a subgroup of $\tilde G$ of index 2 and determine the
structure of the subgroup. We consider the following three cases separately:

Case A:\quad~ $|\wt{\sigma}_p| =|\sigma_p|~\text{ and }~|\wt{\sigma}_q| =|\sigma_q|$;

Case B:\quad~ $|\wt{\sigma}_p|=2|\sigma_p|,~~ |\wt{\sigma}_q| =|\sigma_q|~\text{ or }~
|\wt{\sigma}_p| =|\sigma_p|,~~|\wt{\sigma}_q| =2|\sigma_q|$;

Case C:\quad~ $|\wt{\sigma}_p| =2|\sigma_p|~\text{ and }~ |\wt{\sigma}_q| =2|\sigma_q|$.\\
All the three cases may happen. In fact, the case (A) happens if and only if
$(\frac{p^*}q)=(\frac{q^*}p)=1$; the case (B) happens if and only if
$(\frac{p^*}q)\ne(\frac{q^*}p)$ or $(p,q)=(-1,2)$; and the case (C) happens if and only if
$(\frac{p^*}q)=(\frac{q^*}p)=-1$.

In the case A, we define the subgroup $N$ of $\wt G$ to be
\[
N=\begin{cases}<\wt\sigma_{-1},\wt\sigma_2,\wt\sigma_q^2,\veps> &~~\text{ if }~~p=2\\
<\wt\sigma_p,\wt{\sigma}^2_q,\veps> &~~\text{ if }~~p\ne 2
\end{cases}\eqno(A2.1)
\]
It is easy to see that the subgroup $N$ is abelian of index 2 and is a direct sum of the cyclic
groups generated by the elements. Thus we have
\[
N\cong\begin{cases}
\Z/2\Z\oplus \Z/((q-1)/2)\Z\oplus\Z/2\Z &~~\text{ if }~~p=-1\\
\Z/2\Z\oplus\Z/2\Z\oplus\Z/((q-1)/2)\Z\oplus\Z/2\Z &~~\text{ if }~~p=2\\
\Z/(p-1)\Z\oplus\Z/((q-1)/2)\Z\oplus\Z/2\Z &~~\text{ if }~~p>2.
\end{cases}\eqno (A2.2)
\]

In the case B, we define the the subgroup $N$ of $\wt G$ to be
\[
N=\begin{cases}
<\wt\sigma_{-1},\wt{\sigma}_2,\wt{\sigma}^2_q> &~~\text{ if }~~p=2\\
<\wt{\sigma}_p,\wt{\sigma}^2_q> &~~\text{ if }~~p\ne 2\text{ and }
|\wt{\sigma}_q| =2|\sigma_q|\\
<\wt{\sigma}^2_p,\wt{\sigma}_q> &~~\text{ if }~~|\wt{\sigma}_p| =2|\sigma_p|.
\end{cases}\eqno(B2.1)
\]
Again $N$ is abelian and has index 2 in $\wt G$. In addition, we have
\[
N\cong\begin{cases} \Z/2\Z\oplus\Z/2\Z &~\text{ if }~(p,q)=(-1,2)\\
\Z/2\Z\oplus\Z/(q-1)\Z &~\text{ if }~p=-1,~ q>2\\
\Z/2\Z\oplus\Z/2\Z\oplus\Z/(q-1)\Z &~\text{ if }~p=2\\
\Z/(p-1)\Z\oplus\Z/(q-1)\Z &~\text{ if }~p>2.
\end{cases}\eqno(B2.2)
\]

In the case C, we must have $p,q$ are odd prime numbers. Let $v_2(p-1)$ denote the power
of 2 in $p-1$. We define the subgroup $N$ of $\wt G$ to be
\[
N=\begin{cases}
<\wt{\sigma}^2_p\,,\wt{\sigma}_q> &~~\text{ if }~~v_2(p-1)\le v_2(q-1)\\
<\wt{\sigma}_p\,,\wt{\sigma}^2_q> &~~\text{ if }~~v_2(p-1)> v_2(q-1).
\end{cases}\eqno(C2.1)
\]
Then $N$ is an abelian subgroup of $\wt G$. When $v_2(p-1)\le v_2(q-1)$, we have
\[
|N|=\frac{|\wt\sigma_p^2|\cdot|\wt\sigma_q|}{|<\wt\sigma_p^2>\cap<\wt\sigma_q>|}
=\frac{(p-1)\cdot 2(q-1)}{2}.
\]
So $[\wt G : N]=2$ and $N$ is a normal subgroup of $\wt G$. We have the same result
when $v_2(p-1)> v_2(q-1)$. The subgroup
$<\wt{\sigma}^2_p\,,\wt{\sigma}_q>$ is always an abelian subgroup of $\wt G$ of index 2.
But we can not get all irreducible representations from the inducement
of the representations of this subgroup when $v_2(p-1)> v_2(q-1)$. So we define $N$ in two cases.

Next we determine the structure of the subgroup $N$ in the case C. We consider the
case $v_2(p-1)\le v_2(q-1)$ in detail.
Let $d=\gcd(\frac{p-1}2,q-1)$, $s=(p-1)/2d$ and $t=(q-1)/d$. Choose
$u,v\in\Z$ such that $us+vt=1$. We have the relations
\[(\wt{\sigma}^2_p)^{p-1}=1,~~~
(\wt{\sigma}^2_p)^{\frac{p-1}2}=\veps=\wt{\sigma}^{q-1}_q.
\]

Let $M$ be the free abelian group generated by two words $\alpha\,,\,\beta$. Let
\[
\alpha_1=(p-1)\alpha\;\;;\;\;\beta_1=\frac{p-1}2\alpha-(q-1)\beta\;;
\]
and let $M_1$ be the subgroup of $M$ generated by
$\alpha_1\,,\beta_1$. Then $M_1$ is the kernel of the homomorphism
\[
M\lra N;~~\alpha\mapsto\wt{\sigma}^2_p,~\beta\mapsto\wt{\sigma}_q\;\;.
\]
So we have $N\cong M/M_1$. Define the matrix
\[
A:=\matr{p-1}{\frac{p-1}2}{0}{1-q}\;.
\]
Then $(\alpha_1\,,\beta_1)=(\alpha\,,\beta)\cdot A$.
We determine the structure of $M_1$ by considering the standard form
of $A$. Define
\[
P:=\matr{u}{v}{-t}{s}\in\mathrm{SL}_2(\Z)\;;
\quad\;Q:=\matr{1}{2tv-1}{-1}{-2tv+2}\in\mathrm{SL}_2(\Z).
\]
Then
\[
B:=PAQ=\matr{d}{0}{0}{-2s(q-1)}
\]
is the standard form of $A$. Let
\[
(\tau\,,\mu)=(\alpha\,,\beta)P^{-1}\quad\text{and}\quad(\tau_1\,,\mu_1)=(\alpha_1\,,\beta_1)Q.
\]
Then $(\tau_1\,,\mu_1)=(\tau\,,\mu)B$, $M=\Z\tau\oplus\Z\mu$ and
$M_1=\Z d\tau\oplus \Z 2s(q-1)\mu$. We thus have
\[
N\cong M/M_1\cong \Z/d\Z\oplus\Z/2s(q-1)\Z\;.
\]
By abuse of notation, we also write
\[
(\tau\,,\mu)=(\wt{\sigma}^2_p\,,\wt{\sigma}_q)P^{-1}=
(\wt{\sigma}^{2s}_p\wt{\sigma}^{t}_q\,,\,\wt{\sigma}_p^{-2v}\wt{\sigma}^u_q)\;.
\]
Then $\tau\,,\mu$ are of order $d\,,2s(q-1)$ respectively, and $N$ is a direct sum of
$<\tau>$ and $<\mu>$. We have $\wt\sigma_p^2=\tau^u\mu^{-t}$ and $\wt\sigma_q=\tau^v\mu^s$.
When $v_2(p-1)>v_2(q-1)$, we get the structure of $N$ in the same way. So in the case (C)
we have
\[
N\cong\begin{cases}
\Z/d\Z\oplus\Z/2s(q-1)\Z\ &~~\text{ if }~~v_2(p-1)\le v_2(q-1)\\
\Z/d'\Z\oplus\Z/2s'(p-1)\Z\ &~~\text{ if }~~v_2(p-1)> v_2(q-1),
\end{cases}\eqno(C2.2)
\]
where $d=\gcd(\frac{p-1}2,q-1), s=(p-1)/2d$ and $d'=\gcd(p-1,\frac{q-1}2), s'=(q-1)/2d'$.

Now we summarize our results in the following
\begin{prop}\label{thm2p1}
The abelian subgroup $N$ of the group $\wt G$ of index 2 defined in (A2.1), (B2.1) and (C2.1)
has the structure described in (A2.2), (B2.2) and (C2.2) in the cases (A), (B) and (C), respectively.
In particular, every irreducible representation of $\wt G$ has dimension 1 or 2.
\end{prop}

\section{2-dimensional representations}
We determine all irreducible representations of $\wt G$ in this section. We will use
some basic facts from representation theory freely. For the details, see [6].

It is well-known that the 1-dimensional representations of $\wt G$ correspond bijectively to
those of the maximal abelian quotient $G$ of $\wt G$, which are Dirichlet characters. So we mainly
construct the 2-dimensional irreducible representations of $\wt G$. From the dimension formula of
all irreducible representations, we see that $\wt G$ has $|G|/4$ irreducible representations of
dimension 2, up to isomorphism. Let $N$ be the subgroup of $\wt G$ defined in last section.
Let $\wt G=N\cup\sigma N$ be a decomposition of cosets.
If $\rho:N\rightarrow\Bbb C^*$ is a representation of $N$, the induced representation $\wt\rho$
of $\rho$ is a representation of $\wt G$ of dimension $2$. The space of the representation
$\wt\rho$ is $V=\text{Ind}_N^{\wt G}(\Bbb C)=\Bbb C[\wt G]\otimes_{\Bbb C[N]}\Bbb C$ with basis
$e_1=1\otimes 1$ and $e_2=\sigma\otimes 1$. The group homomorphism
\[
\wt\rho : \wt G\lra \GL(V)\simeq\GL_2(\C)
\]
is given under the basis by
\begin{equation}\label{eq2p2}
\wt\rho(\wt\sigma)=\matr{\rho(\wt\sigma)}{\rho(\wt\sigma\sigma)}
{\rho(\sigma^{-1}\wt{\sigma})}{\rho(\sigma^{-1}\wt{\sigma}\sigma)},\quad\;\forall\;\wt\sigma\in\wt
G,
\end{equation}
where $\rho(\wt\sigma)=0$ if $\wt\sigma\notin N$. The representation $\wt\rho$ is irreducible
if and only if $\rho\not\cong \rho^{\tau}$ for every $\tau\in\wt G\setminus N$,
where $\rho^{\tau}$ is the conjugate representation of $\rho$ defined by
\[
\rho^{\tau}(x)=\rho(\tau^{-1}x\tau)\;,\;\;\forall\;x\in N\,.
\]
Since $N$ is abelian, we only need to check $\rho\not\cong \rho^{\sigma}$.

Now we begin to construct all 2-dimensional irreducible representations of $\wt G$.
As in last section,
we consider the three cases separately. In addition, we consider the case when $p$ and $q$
are odd prime numbers in details, and only state the results in the cases when $p=-1$ or $2$.

\noindent {\bf 3.1. Case A.} Assume $p>2$. We have in this case
$N=\langle\wt{\sigma}_p\,,\,\wt{\sigma}^2_q\,,\veps\rangle$
and
\[
N\cong\Z/(p-1)\Z\oplus\Z/((q-1)/2)\Z\oplus\Z/2\Z\;.
\]
Every irreducible representation of $N$ can be written as
$\rho_{ijk}: N\lra\C^*$ with
\[
\rho_{ijk}(\wt{\sigma}_p)=\z^i_{p-1}\;;\;
\rho_{ijk}(\wt{\sigma}^2_q)=\z^{2j}_{q-1}\;;\;
\rho_{ijk}(\veps)=(-1)^k\;.
\]
where $0\le i<p-1,~0\le j<\frac{q-1}2$ and $k=0,1$.
Since $\wt G=N\cup\wt\sigma_q N$ and $\rho^{\wt{\sigma}_q}_{ijk}(\wt{\sigma}_p)
=\rho_{ijk}(\veps)\rho_{ijk}(\wt{\sigma}_p)=(-1)^k\rho_{ijk}(\wt{\sigma}_p)$, we have
\[
\rho^{\wt{\sigma}_q}_{ijk}\not\cong\rho_{ijk}\Longleftrightarrow
k=1.
\]
Write $\rho_{ij}=\rho_{ij1}$. The induced representation $\wt\rho_{ij}:
\wt{G}\longrightarrow \mathrm{GL}_2(\C)$ of $\rho_{ij}$ is given by
\[
\wt\rho_{ij}(\wt{\sigma}_p)=\matr{\z_{p-1}^i}{0}{0}{-\z_{p-1}^i},~
\wt\rho_{ij}(\wt{\sigma}_q)=\matr{0}{\z^{2j}_{{q-1}}}{1}{0},~
\wt\rho_{ij}(\veps)=-I,\eqno(A3.1)
\]
where $I$ is the identity matrix of degree 2. Since
\[
\wt\rho_{ij}(\wt{\sigma}^2_p)=\matr{\z^{2i}_{p-1}}{0}{0}{\z^{2i}_{p-1}}\quad\text{and}\quad
\wt\rho_{ij}(\wt{\sigma}^2_q)=\matr{\z^{2j}_{{q-1}}}{0}{0}{\z^{2j}_{{q-1}}},
\]
we see that the representations $\wt\rho_{ij}$ with $0\le i<\frac{p-1}2,~0\le j<\frac{q-1}2$
are irreducible and are not isomorphic to each other, by considering the values of the characters
of these representations at $\wt\sigma_p^2$ and $\wt\sigma_q^2$. The number of these representations
is $\frac{p-1}{2}\cdot\frac{q-1}{2}=\frac{|G|}{4}$. So they are all the irreducible
representations of $\wt G$ of dimension 2.

Similarly, when $p=-1$, all irreducible representations of $\wt G$ of dimension 2 are
$\wt\rho_{j}$ with $0\le j<\frac{q-1}2$, where
\[
\wt\rho_j(\wt\sigma_{-1})=\matr{1}{0}{0}{-1},\quad
\wt\rho_{j}(\wt{\sigma}_q)=\matr{0}{\z^{2j}_{{q-1}}}{1}{0},\quad
\wt\rho(\veps)=-I\eqno(A3.2)
\]
and when $p=2$, all irreducible representations of $\wt G$ of dimension 2 are
$\bar\rho_{ij}$ with $0\le i\le 1$ and $0\le j<\frac{q-1}2$, where $\bar\rho_{ij}(\veps)=-I$
and
\[
\bar\rho_{ij}(\wt\sigma_{-1})=(-1)^iI,~\bar\rho_{ij}(\wt\sigma_{2})=\matr{1}{0}{0}{-1},~
\bar\rho_{ij}(\wt{\sigma}_q)=\matr{0}{\z^{2j}_{{q-1}}}{1}{0}.
\eqno (A3.3)\]

\noindent {\bf 3.2 Case B.} Assume $p>2$ and $|\wt\sigma_q|=2|\sigma_q|$. Then
$N=\langle\wt{\sigma}_p\,,\,\wt{\sigma}^2_q\rangle$, and
\[
N\cong\Z/(p-1)\Z\oplus\Z/(q-1)\Z.
\]
Any irreducible representation of $N$ has the form $\rho_{ij}: N\lra\C^*$, where
\[
\rho_{ij}(\wt{\sigma}_p)=\z^i_{p-1}\;,\;
\rho_{ij}(\wt{\sigma}^2_q)=\z^j_{q-1}\;,\;
\rho_{ij}(\veps)=\rho_{ij}(\wt{\sigma}^2_q)^{\frac{q-1}2}=(-1)^j\;,
\]
and $0\le i<p-1,~0\le j<q-1$. It is easy to check that
\[
\rho^{\wt{\sigma}_q}_{ij}\not\cong\rho_{ij}\Longleftrightarrow
j\equiv 1\;\pmod{2}\;.
\]
The induced representation $\wt\rho_{ij}:
\wt{G}\longrightarrow \mathrm{GL}_2(\C)$ of $\rho_{ij}$ with odd $j$ is given by
\[
\wt\rho_{ij}(\wt{\sigma}_p)=\matr{\z^i_{p-1}}{0}{0}{-\z^i_{p-1}},~~
\wt\rho_{ij}(\wt{\sigma}_q)=\matr{0}{\z^j_{q-1}}{1}{0}.\eqno(B3.1)
\]
Since
\[
\wt\rho_{ij}(\wt{\sigma}^2_p)=\matr{\z^{2i}_{p-1}}{0}{0}{\z^{2i}_{p-1}}\quad\text{and}\quad
\wt\rho_{ij}(\wt{\sigma}^2_q)=\matr{\z^{j}_{q-1}}{0}{0}{\z^{j}_{q-1}},
\]
we see that the representations $\wt\rho_{ij}$ with $0\le i<\frac{p-1}2$ and $0\le j< q-1,~2\nmid j$
are irreducible and are not isomorphic to each other. The number of these representations is
$\frac{|G|}4$. So they are all the irreducible representations of $\wt G$ of dimension 2.

Similarly, when $(p,q)=(-1,2)$, there is only one irreducible representation $\wt\rho_0$ of dimension 2
defined by
\[
\wt\rho_0(\wt\sigma_{-1})=\matr{1}{0}{0}{-1},\quad\text{and}\quad
\wt\rho_0(\wt\sigma_{2})=\matr{0}{-1}{1}{0}.\eqno(B3.2)
\]

When $p=-1$ and $q>2$, all irreducible representations of dimension 2 are $\wt\rho_{j}$ with
$0\le j< q-1,~2\nmid j$, where $\wt\rho_j$ is defined by
\[
\wt\rho_j(\wt\sigma_{-1})=\matr{1}{0}{0}{-1}\quad\text{and}\quad
\wt\rho_{j}(\wt{\sigma}_q)=\matr{0}{\z^j_{q-1}}{1}{0}. \eqno(B3.3)
\]

When $p=2$, all irreducible representations of dimension 2 are $\bar\rho_{ij}$ with
$0\le i\le 1$ and $0\le j< q-1,~2\nmid j$, where $\bar\rho_{ij}$ is defined by
\[
\bar\rho_{ij}(\wt\sigma_{-1})=(-1)^iI,~
\bar\rho_{ij}(\wt\sigma_{2})=\matr{1}{0}{0}{-1},~
\bar\rho_{ij}(\wt{\sigma}_q)=\matr{0}{\z^j_{q-1}}{1}{0}.
\eqno(B3.4)
\]

When $|\wt\sigma_p|=2|\sigma_p|$, all irreducible representations of dimension 2 are $\hat\rho_{ij}$
with $0\le i< p-1,~2\nmid i$ and $0\le j<\frac{q-1}2$, where $\hat\rho_{ij}$ is defined by
\[
\hat\rho_{ij}(\wt{\sigma}_p)=\matr{0}{\z^i_{p-1}}{1}{0},~~
\hat\rho_{ij}(\wt{\sigma}_q)=\matr{\z^j_{q-1}}{0}{0}{-\z^j_{q-1}}.\eqno(B3.5)
\]

\noindent {\bf 3.3. Case C.} Assume $v_2(p-1)\le v_2(q-1)$. Let
\[
d=\gcd(\frac{p-1}{2}\,,q-1)\;,\;s=\frac{p-1}{2d}\,,t=\frac{q-1}d\;;\;\;us+vt=1
\]
as before. We must have that $t$ is even and $u$ is odd. Let
$\tau=\wt{\sigma}^{2s}_p\cdot\wt{\sigma}^{t}_q$
and $\mu=\wt{\sigma}^{-2v}_p\cdot\wt{\sigma}^u_q$.
Then
$N=\langle\wt{\sigma}^2_p\,,\,\wt{\sigma}_q\rangle=\langle\tau\,,\mu\rangle$
and
\[
N\cong\Z/d\Z\oplus\Z/2s(q-1)\Z.
\]
Any irreducible representation $\rho_{ij}: N\lra\C^*$ is of the form
\[
\rho_{ij}(\tau)=\z^i_d=\z^{2s(q-1)i}_{(p-1)(q-1)}\quad\text{and}\quad
\rho_{ij}(\mu)=\z^j_{2s(q-1)}=\z^{dj}_{(p-1)(q-1)}.
\]
From $\wt\sigma_p^2=\tau^u\mu^{-t}$ and $\wt\sigma_q=\tau^v\mu^s$, we have
\[
\rho_{ij}(\wt{\sigma}^2_p)=\z^{2sui-j}_{p-1}\;;\;
\rho_{ij}(\wt{\sigma}_q)=\z^{2tvi+j}_{2(q-1)}\;;\;
\rho_{ij}(\veps)=\rho_{ij}(\wt{\sigma}^2_p)^{\frac{p-1}2}=(-1)^j\;.
\]
It is easy to show
\[
\rho^{\wt{\sigma}_p}_{ij}\not\cong\rho_{ij}\Longleftrightarrow
j\equiv 1\;\pmod{2}\;.
\]
The induced representation $\wt\rho_{ij}:\wt{G}\longrightarrow \mathrm{GL}_2(\C)$
of $\rho_{ij}$ with odd $j$ is given by
\[
\wt\rho_{ij}(\tau)=\matr{\z^i_d}{0}{0}{\z^i_d}\;;\quad\;
\wt\rho_{ij}(\mu)=\matr{\z^j_{2s(q-1)}}{0}{0}{-\z^j_{2s(q-1)}}.
\]
Here in the first equality we used the fact that $t$ is even, and in the second equality
we used the fact that $u$ is odd. Furthermore we have
\[
\wt\rho_{ij}(\wt{\sigma}_p)=\matr{0}{\z^{2sui-j}_{p-1}}{1}{0}\;;\;
\wt\rho_{ij}(\wt{\sigma}_q)=\matr{\z^{2tvi+j}_{2(q-1)}}{0}{0}{-\z^{2tvi+j}_{2(q-1)}}.
\eqno(C3.1)\]
By considering the values of the character of $\wt\rho_{ij}$ at $\tau$ and $\mu^2$, we see that
all the representations $\wt\rho_{ij}$ with $0\le i<d$ and $0\le j<s(q-1),~2\nmid j$ are
irreducible and are not isomorphic to each other. The number of these representations is
$d\cdot \frac{s(q-1)}{2}=\frac{|G|}{4}$. So they are all the irreducible representations
of $\wt G$ of dimension 2.

Similarly, if $v_2(p-1)> v_2(q-1)$, we let
\[
d'=\gcd(p-1\,,\frac{q-1}2)\;,\;s'=\frac{p-1}{d}\,,t'=\frac{q-1}{2d}\;;\;\;u's'+v't'=1.
\]
Then all the irreducible representations of $\wt G$ of dimension 2 are $\hat\rho_{ij}$ with
$0\le i<d'$ and $0\le j<t'(p-1),~2\nmid j$, where $\hat\rho_{ij}$ is defined by
\[
\hat\rho_{ij}(\wt{\sigma}_p)=\matr{\z^{2s'u'i+j}_{2(p-1)}}{0}{0}{-\z^{2s'u'i+j}_{2(p-1)}};~~
\hat\rho_{ij}(\wt{\sigma}_q)=\matr{0}{\z^{2t'v'i-j}_{q-1}}{1}{0}\; .
\eqno(C3.2)
\]

Let R$^2(\wt G)$ be the set of all irreducible representations, up to isomorphism, of
$\wt G$ of dimension 2. As a summary, we have proved the following.
\begin{thm} All 2-dimensional irreducible representations of $\wt G$ are induced from the
representations of $N$. In detail, we have

In the case (A)
\[
\mathrm{R}^2(\wt G)=\begin{cases}
\{\wt\rho_j\mid 0\le j<\frac{q-1}2\}&\mathrm{ if }~p=-1\\
\{\bar\rho_{ij}\mid i=0,1,~0\le j<\frac{q-1}2\} &\mathrm{ if }~p=2\\
\{\wt\rho_{ij}\mid 0\le i<\frac{p-1}2,~0\le j<\frac{q-1}2\} &\mathrm{ if }~p>2,
\end{cases}
\]
where $\wt\rho_j,~\bar\rho_{ij}$ and $\wt\rho_{ij}$ are defined in (A3.2), (A3.3) and
(A3.1) respectively.

In the case (B)
\[
\mathrm{R}^2(\wt G)=\begin{cases}
\{\wt\rho_0\}&\mathrm{ if }~(p,q)=(-1,2)\\
\{\wt\rho_j\mid 0\le j<q-1,~2\nmid j\}&\mathrm{ if }~p=-1,~q>2\\
\{\bar\rho_{ij}\mid i=0,1,~0\le j<q-1,~2\nmid j\} &\mathrm{ if }~p=2\\
\{\hat\rho_{ij}\mid 0\le i<p-1,~2\nmid i,~0\le j<\frac{q-1}2\} &\mathrm{ if }~|\wt\sigma_p|=2|\sigma_p|\\
\{\wt\rho_{ij}\mid 0\le i<\frac{p-1}2,~0\le j<q-1,~2\nmid j\} &\mathrm{ otherwise },
\end{cases}
\]
where $\wt\rho_0,~\wt\rho_j,~\bar\rho_{ij},~\hat\rho_{ij}$  and $\wt\rho_{ij}$ are defined in
(B3.2), (B3.3), (B3.4), (B3.5) and (B3.1) respectively.

In the case (C)
\[
\mathrm{R}^2(\wt G)=\begin{cases}
\{\wt\rho_{ij}\mid 0\le i<d,~0\le j<s(q-1),~2\nmid j\} &\mathrm{ if }~v_2(p-1)\le v_2(q-1),\\
\{\hat\rho_{ij}\mid 0\le i<d',~0\le j<t'(p-1),~2\nmid j\} &\mathrm{ otherwise},
\end{cases}
\]
where $\wt\rho_{ij}$ and $\hat\rho_{ij}$ are defined in (C3.1) and (C3.2) respectively.
\end{thm}

\section{The Frobenious maps}
This section is a preparation for the next section to compute the Artin $L$-functions of the
quasi-cyclotomic fields $\wt K$ when $p=-1$. For a prime number $\ell$ which is unramified
in $\wt K/K$, let $I_\ell$ (resp. $\wt I_\ell$) be the inertia group of $\ell$ in the extension
$K/\Bbb Q$ (resp. $\wt K/\Bbb Q$). Let ${\mathrm{Fr}}_\ell$ be the Frobenious automorphism of $\ell$
in $G/I_\ell$ and $\wt{\mathrm{Fr}}_\ell$ the Frobenious automorphism of $\ell$ in $\wt G/\wt I_\ell$
associated to some prime ideal over $\ell$. In this section we determine $\wt{\mathrm{Fr}}_\ell$ by
Fr$_\ell$ for $\ell=2$.

From now on, we always assume that $p=-1$, namely, $K=\Bbb Q(\zeta_{4q})$ and $\wt K=K(\sqrt[4]{q^*})$.
For a prime number $\ell$, we say that $\ell$ is ramified (resp. inertia, splitting)
in the relative quadratic extension $\tilde K/K$
if the prime ideals of $K$ over $\ell$ are ramified (resp. inertia, splitting) in $\tilde K$.
In [Sect.5, 7] we have determined the decomposition nature of odd prime numbers in $\wt K/K$. Now we determine
the decomposition nature of 2 in $\wt K/K$.
\begin{prop} If $q=2$, then 2 is ramified in $\wt K/K$. If $q$ is odd, then 2 is unramified in $\wt K/K$
if and only if $(\frac 2q)=1$, and in this case 2 splits in $\wt K/K$ if $q^*\equiv 1\mod 16$ and is inertia
in $\wt K/K$ if $q^*\equiv 1\mod 8$ but $q^*\not\equiv 1\mod 16$.
\end{prop}
\begin{proof} We first consider the case $q=2$. The unique prime ideal of $K$ over 2 is the principal
ideal generated by $\pi:=1-\zeta_8$. Since the ramification degree of 2 in $K/\Bbb Q$ is 4 and
$\sqrt 2=\pi(\pi+2\zeta_8)\zeta_8$, we have that 2 is ramified in $\wt K/K$ if and only if
$x^2\equiv\sqrt 2\mod\pi^{10}$ is not solvable in the ring $O_K$ of the integers of $K$ by [2], which
is equivalent to that $(1+\frac 2\pi\zeta_8)\zeta_8\mod\pi^{8}$ is not a square. Since $2=u\pi^4$ for some
unit $u$, we have
\[
(1+\frac 2\pi\zeta_8)\zeta_8\equiv\zeta_8\equiv (1-\pi)\mod\pi^3,
\]
namely $(1+\frac 2\pi\zeta_8)\zeta_8\mod\pi^3$ is not a square. So 2 is ramified in $\wt K/K$.

Now we assume that $q$ is odd. Let $\pi_2=1-\zeta_4$. Then 2 is unramified in $\wt K/K$ if and only
if $x^2\equiv\sqrt{q^*}\mod\pi_2^4$ is solvable in $O_K$. Furthermore, 2 splits in $\wt K/K$ if and
only if $x^2\equiv\sqrt{q^*}\mod\pi_2^5$ is solvable in $O_K$. By Gauss sum we have
\[
\sqrt{q^*}=\sum_{a=1}^{q-1}(\frac aq)\zeta_q^a=1+2\sum_{(\frac aq)=1}\zeta_q^a.
\]

Let $\alpha=\sum_{(\frac aq)=1}\zeta_q^a,~\beta=\sum_{(\frac aq)=1}\zeta_{2q}^a$, and
$\gamma=\sum_{(\frac aq)=1}\sum_{(\frac bq)=1,~a<b}\zeta_{2q}^{a+b}$,
where in the summations $a,b$ run over $1,2,\cdots,q-1$.
Then $\alpha=\beta^2-2\gamma$, from which and the equality $2=\pi_2^2-\pi_2^3$, we have
\[\begin{split}
\sqrt{q^*}&=1+2\beta^2-4\gamma=1+\pi_2^2\beta^2-\pi_2^3\beta^2-4\gamma\\
&\equiv(1+\pi_2\beta)^2-\pi_2^3(\beta+\beta^2)+\pi_2^4(\beta-\gamma)\\
&\equiv(1+\pi_2\beta)^2-\pi_2^3(\alpha+\beta)+\pi_2^4(\beta+\gamma)\mod\pi_2^5.
\end{split}\]
Since $\zeta_{2q}=-\zeta_q^{-\frac{q-1}2}=-\zeta_q^t$, where $t$ is the inverse of 2 in $(\Z/q\Z)^*$,
we see $\beta=\sum_{(\frac aq)=1}(-1)^a\zeta_{q}^{ta}\equiv\sum_{(\frac aq)=1}\zeta_{q}^{ta}\mod 2$.
So if $(\frac 2q)=1$ we have $\alpha\equiv\beta\mod 2$ and thus $2$ is unramified in $\wt K/K$,
and if $(\frac 2q)=-1$ we have $\alpha+\beta\equiv\sum_{a=1}^{q-1}\zeta_q^a=-1\mod 2$ and thus
$2$ is ramified in $\wt K/K$.

Now we assume $(\frac 2q)=1$. Then $\sqrt{q^*}\mod\pi_2^5$ is a square if and only if
$\pi_2\mid\beta+\gamma$. We consider $2(\beta+\gamma)$. Since $\alpha\equiv\beta\mod 2$, we have
\[
2(\beta+\gamma)=2\beta+\beta^2-\alpha\equiv \alpha(\alpha+1)\mod 4.
\]
From $\sqrt{q^*}=1+2\alpha$, we see $\alpha(\alpha+1)=\frac{q^*-1}4$. Since $8\mid q^*-1$ under
the assumption $(\frac 2q)=1$, we have $\beta+\gamma\equiv\frac{q^*-1}8\mod 2$. So
$\pi_2\mid\beta+\gamma$ if and only if $\pi_2\mid\frac{q^*-1}8$, namely $2\mid\frac{q^*-1}8$.
We complete the proof.
\end{proof}

By the way, we have determined the ring $O_{\wt K}$ of the integers of $\wt K$. In fact, we have
\begin{cor}
Assume that $q$ is an odd prime number. Let $t_q=(q-1)/4$ if $q\equiv 1\mod 4$, and
$t_q=(q-3)/4$ if $q\equiv 3\mod 4$. Then
\[
O_{\wt K}=\begin{cases}
\Bbb Z\left[\zeta_{4q},\frac{\sqrt[4]{q^*}+1+\pi_2\beta}{\pi_2(\sin[\frac 1q])^{t_q}}\right]
&\mathrm{if}~ (\frac 2q)=-1\\
\Bbb Z\left[\zeta_{4q},\frac{\sqrt[4]{q^*}+1+\pi_2\beta}{2(\sin[\frac 1q])^{t_q}}\right]
&\mathrm{if}~ (\frac 2q)=1.
\end{cases}\]
\end{cor}
\begin{proof}
See [Th.2, 7].
\end{proof}

Now we assume that 2 is unramified in $\wt K/K$. Let Fr$_2\in G$ such
that Fr$_2(\zeta_4)=1$ and Fr$_2(\zeta_q)=\zeta_q^2$. It is a Frobenious element of 2 in
$G$ modulo $I_2$. We have Fr$_2=\sigma_2^{b_2}$ with $2\mid b_2$ for $(\frac 2q)=1$. Thus
$\wt{\text{Fr}}_2=\wt\sigma_2^{b_2}$ or $\wt{\text{Fr}}_2=\wt\sigma_2^{b_2}\varepsilon$.
We need to determine $\wt{\text{Fr}}_2$ completely. Since $(\frac 2q)=1$, we have
\[
\sqrt{q^*}\equiv(1+\pi_2\alpha)^2+\pi_2^4(\beta+\gamma)\mod\pi_2^5.
\]
Write $u=1+\pi_2\alpha$ for simplicity. Since $\sqrt{q^*}\equiv u^2\mod\pi_2^4$, we see
$\frac{\sqrt[4]{q^*}-u}2\in O_{\wt K}$. Let $\wp$ be the prime ideal of $\wt K$ over 2
associated to $\wt{\text{Fr}}_2$. By the definition, we have
\[
\wt{\text{Fr}}_2\left(\frac{\sqrt[4]{q^*}-u}2\right)\equiv\left(\frac{\sqrt[4]{q^*}-u}2\right)^2
\equiv(\beta+\gamma)+\frac{\sqrt[4]{q^*}-u}2\mod\wp.
\]
On the other hand, since $\wt\sigma_q^{b_2}(\sqrt[4]{q^*})=(-1)^{\frac{b_2}2}\sqrt[4]{q^*}$ and
$\wt\sigma_q^{b_2}(u)=u$ for $2\mid b_2$, we have
\[
\wt\sigma_q^{b_2}(\frac{\sqrt[4]{q^*}-u}2)=\frac{(-1)^{\frac{b_2}2}\sqrt[4]{q^*}-u}2
\]
and
\[
\wt\sigma_q^{b_2}\varepsilon(\frac{\sqrt[4]{q^*}-u}2)=\frac{(-1)^{\frac{b_2}2+1}\sqrt[4]{q^*}-u}2.
\]
So if $2\mid\frac{b_2}2$ we have $\wt{\text{Fr}}_2=\wt\sigma_2^{b_2}$ if and only if
$\pi_2\mid\beta+\gamma$ (namely 2 splits in $\wt K/K$), and if $2\nmid\frac{b_2}2$ we have
$\wt{\text{Fr}}_2=\wt\sigma_2^{b_2}$ if and only if $\pi_2\nmid\beta+\gamma$ (namely 2 is inertia
in $\wt K/K$). In the case $q\equiv 3\mod 4$, we can always assume that $2\nmid\frac{b_2}2$,
since if $4\mid b_2$, we may replace $b_2$ by $b_2+(q-1)$. In the case $q\equiv 1\mod 4$, we have
$2\mid\frac{b_2}2\Longleftrightarrow 2^{\frac{q-1}4}\equiv 1\mod q\Longleftrightarrow q$ has the
form $A^2+64B^2$ for $A,B\in\Bbb Z$, by the Exercise 28 in Chap.5 in [5]. So we get the following result
\begin{prop} Assume that 2 is unramified in $\wt K/K$. Let $\mathrm{Fr}_2=\sigma_2^{b_2}$. We have
$2\mid b_2$. If $q\equiv 3\mod 4$, we always assume $b_2\equiv 2\mod 4$. Let $P_0$ be the set of
the prime numbers of the form $A^2+64B^2$ with $A,B\in\Bbb Z$. Then we have
\[
\wt{\mathrm{Fr}}_2=\begin{cases}\wt\sigma_2^{b_2}&\mathrm{if}~~q\not\in P_0,16\nmid q^*-1,
~\mathrm{or}~q\in P_0,16\mid q^*-1\\
\wt\sigma_2^{b_2}\varepsilon &\mathrm{if}~~q\in P_0,16\nmid q^*-1,
~\mathrm{or}~q\not\in P_0,16\mid q^*-1.
\end{cases}
\]
\end{prop}

The following lemma is useful in next section.
\begin{lemma} We have $\varepsilon\in\wt I_\ell$ if and only if $\ell$ is ramified in
$\wt K/K$.
\end{lemma}
\begin{proof}
The canonical projection $\wt G\lra G\simeq\wt G/\langle\veps\rangle$
induces a surjective homomorphism $\wt{I}_\ell\longrightarrow I_\ell$ which implies
the isomorphism $\wt{I}_\ell/<\varepsilon>\cap\wt{I}_\ell\cong I_\ell$. Thus
$\ell$ is ramified in $\wt K/K\Longleftrightarrow |\wt{I}_\ell|=2|I_\ell|
\Longleftrightarrow |\wt{I}_\ell\cap<\varepsilon>|=2\Longleftrightarrow\varepsilon\in\wt{I}_\ell$.
\end{proof}
\section{The Artin $L$-functions}
In this section we compute the Artin $L$-functions of the quasi-cyclotomic fields
$\wt K=\Bbb Q(\zeta_{4q},\sqrt[4]{q^*})$.

The $L$-functions associated to the 1-dimensional representations of $\wt G$ are the well-known
Dirichlet $L$-functions. So we mainly compute the $L$-functions associated to the 2-dimensional
irreducible representation of $\wt G$. Let $\varphi:\wt G\rightarrow\text{GL}(V)$ be a 2-dimensional
irreducible representations. The Artin $L$-function $L(\varphi,s)$ associated to $\varphi$ is
defined as the product of the local factors
\[
L(\varphi,s)=\prod_{\ell:\mathrm{prime}}L_\ell(\varphi,s),
\]
where the local factors are defined as
$L_\ell(\varphi,s)=\mathrm{det}(1-\varphi(\wt{\mathrm{Fr}}_\ell)\ell^{-s}|V^{\wt I_\ell})^{-1}$.
Now we begin to compute them. First we notice that if $\ell$ is ramified in $\wt K/K$, then
$V^{\wt I_\ell}=0$ and $L_\ell(\varphi,s)=1$, which is due to the facts that
$\varepsilon\in\wt I_\ell$ by Lem.4.4 and $\varphi(\varepsilon)=-I$ for any irreducible representation
$\varphi$ of $\wt G$ by Th.3.1.

\noindent {\bf 5.1 The case $q=2$.} By section 3, there is only one 2-dimensional representation
$\wt\rho_0$ in this case, which is defined as
\[
\wt\rho_0(\wt\sigma_{-1})=\matr{1}{0}{0}{-1},\quad\text{and}\quad
\wt\rho_0(\wt\sigma_{2})=\matr{0}{-1}{1}{0}.
\]
Since 2 is ramified in $\wt K/K$, we have $L_2(\wt\rho_0,s)=1$. Assume that $\ell$ is an odd prime
number.

If $\ell\equiv 7\mod 8$, then Fr$_\ell=\sigma_{-1}$ and thus $\wt{\text{Fr}}_\ell=\wt\sigma_{-1}$
or $\wt\sigma_{-1}\varepsilon$. In any case we have
\[
L_\ell(\wt\rho_0,s)=\mathrm{det}\left(I\pm\matr{1}{0}{0}{-1}\ell^{-s}\right)^{-1}
=(1-\ell^{-2s})^{-1}.
\]

If $\ell\equiv 5\mod 8$, then Fr$_\ell=\sigma_{2}$ and thus $\wt{\text{Fr}}_\ell=\wt\sigma_{2}$
or $\wt\sigma_{2}\varepsilon$. We have
\[
L_\ell(\wt\rho_0,s)=\mathrm{det}\left(I\pm\matr{0}{-1}{1}{0}\ell^{-s}\right)^{-1}
=(1+\ell^{-2s})^{-1}.
\]

If $\ell\equiv 3\mod 8$, then Fr$_\ell=\sigma_{-1}\sigma_{2}$ and thus
$\wt{\text{Fr}}_\ell=\wt\sigma_{-1}\wt\sigma_{2}$
or $\wt\sigma_{-1}\wt\sigma_{2}\varepsilon$. We have
\[
L_\ell(\wt\rho_0,s)=\mathrm{det}\left(I\pm\matr{1}{0}{0}{-1}\matr{0}{-1}{1}{0}\ell^{-s}\right)^{-1}
=(1-\ell^{-2s})^{-1}.
\]

If $\ell\equiv 1\mod 8$, then Fr$_\ell=1$ and thus $\wt{\text{Fr}}_\ell=1$ or $\varepsilon$.
In this case, we must determine $\wt{\text{Fr}}_\ell$ completely.
Since $\wt{\text{Fr}}_\ell(\sqrt[4]2)\equiv(\sqrt[4]2)^\ell\mod\wp$ for the prime ideal $\wp$
of $\wt K$ over $\ell$ associated to $\wt{\text{Fr}}_\ell$, we have $\wt{\text{Fr}}_\ell=1$
if $2^{\frac{\ell-1}4}\equiv 1\mod\ell$, and $\wt{\text{Fr}}_\ell=\varepsilon$ if
$2^{\frac{\ell-1}4}\equiv -1\mod\ell$. As in last section, we have that for $\ell\equiv 1\mod 8$,
$2^{\frac{\ell-1}4}\equiv 1\mod\ell$ if and only if $\ell\in P_0$. So we have
\[
L_\ell(\wt\rho_0,s)=\begin{cases}
(1-\ell^{-s})^{-2}\quad &\text{if }\ell\in P_0\\
(1+\ell^{-s})^{-2}\quad &\text{otherwise. }
\end{cases}
\]

We get the Artin $L$-function in the case $(p,q)=(-1,2)$ as follows.
\begin{equation}\label{eq2p3}\begin{split}
L(\wt\rho_0,s)=&\prod_{\ell\equiv 3\text{ or }7(8)}(1-\ell^{-2s})^{-1}\cdot
\prod_{\ell\equiv 5(8)}(1+\ell^{-2s})^{-1}\\
&\times\prod_{\ell\in P_0}(1-\ell^{-s})^{-2}\cdot
\prod_{\ell\equiv 1(8),~\ell\not\in P_0}(1+\ell^{-s})^{-2}.
\end{split}\end{equation}

\noindent {\bf 5.2 The case $q$ is odd.} All 2-dimensional irreducible representations
of $\wt G$ are $\wt\rho_j$ with $0\le j<q-1,2\mid j$ if $q\equiv 1\mod 4$, and
$0\le j<q-1,2\nmid j$ if $q\equiv 3\mod 4$, where $\wt\rho_j$ is defined by
\[
\wt\rho_j(\wt\sigma_{-1})=\matr{1}{0}{0}{-1},\quad\wt\rho_{j}(\wt{\sigma}_q)
=\matr{0}{\z^{j}_{q-1}}{1}{0}\quad\text{and}\quad\wt\rho_{j}(\varepsilon)=-I.
\]
We first determine the local factors $L_\ell(\wt\rho_j,s)$ for $\ell\neq 2,q$. For such
$\ell$, we have $V^{\wt I_\ell}=V$. Let Fr$_\ell=\sigma_{-1}^{a_\ell}\sigma_q^{b_\ell}$,
which is equivalent to
$\ell\equiv(-1)^{a_\ell}\mod 4$ and $\ell\equiv g^{b_\ell}\mod q$, where $g$ is the primitive
root $\mod q$ associated to $\sigma _q$. It is easy to compute that
\[
\wt\rho_j(\wt\sigma_q^{b_\ell})=\matr{0}{\zeta_{q-1}^{j}}{1}{0}^{b_\ell}=\begin{cases}
\zeta_{2(q-1)}^{jb_\ell}I&\quad\text{ if }2\mid b_\ell\\
\matr{0}{\zeta_{2(q-1)}^{j(b_\ell+1)}}{\zeta_{2(q-1)}^{j(b_\ell-1)}}{0}
&\quad\text{ if }2\nmid b_\ell.
\end{cases}
\]
Furthermore, we have
\[
\text{det}(I-\wt\rho_j(\wt\sigma_{-1}^{a_\ell}\wt\sigma_q^{b_\ell})\ell^{-s})=\begin{cases}
(1-\zeta_{2(q-1)}^{jb_\ell}\ell^{-s})^2 &\text{ if }a_\ell=0,~2\mid b_\ell\\
1-\zeta_{q-1}^{jb_\ell}\ell^{-2s} &\text{ if }a_\ell=0,~2\nmid b_\ell\\
&\text{ or }a_\ell=1,~2\mid b_\ell\\
1+\zeta_{q-1}^{jb_\ell}\ell^{-2s} &\text{ if }a_\ell=1,~2\nmid b_\ell
\end{cases}
\]
and
\[
\text{det}(I+\wt\rho_j(\wt\sigma_{-1}^{a_\ell}\wt\sigma_q^{b_\ell})\ell^{-s})=\begin{cases}
(1+\zeta_{2(q-1)}^{jb_\ell}\ell^{-s})^2 &\text{ if }a_\ell=0,~2\mid b_\ell\\
1-\zeta_{q-1}^{jb_\ell}\ell^{-2s} &\text{ if }a_\ell=0,~2\nmid b_\ell\\
&\text{ or }a_\ell=1,~2\mid b_\ell\\
1+\zeta_{q-1}^{jb_\ell}\ell^{-2s} &\text{ if }a_\ell=1,~2\nmid b_\ell.
\end{cases}
\]
So we get
\[
L_\ell(\wt\rho_j,s)=(1-\zeta_{q-1}^{jb_\ell}\ell^{-2s})^{-1}
\]
if $\ell\equiv 1\mod 4$ and $\ell\equiv g^{b_\ell}\mod q$ with $2\nmid b_\ell$, or if
$\ell\equiv 3\mod 4$ and $\ell\equiv g^{b_\ell}\mod q$ with $2\mid b_\ell$, and
\[
L_\ell(\wt\rho_j,s)=(1+\zeta_{q-1}^{jb_\ell}\ell^{-2s})^{-1}
\]
if $\ell\equiv 3\mod 4$ and $\ell\equiv g^{b_\ell}\mod q$ with $2\nmid b_\ell$.

To compute the local factors when $\ell\equiv 1\mod 4$ and $\ell\equiv g^{b_\ell}\mod q$
with $2\mid b_\ell$, we must determine $\wt{\text{Fr}}_\ell$ completely. Since $(\frac\ell q)=1$,
we have $(\frac q\ell)=1$ and $(\frac{q^*}\ell)=1$. Let $\alpha_\ell\in\Z$ such that
$\alpha_\ell^2\equiv q^*\mod\ell$. From
$\wt\sigma_q^{b_\ell}(\sqrt[4]{q^*})=(-1)^{\frac{b_\ell}2}\sqrt[4]{q^*}$, we see
$\wt{\text{Fr}}_\ell=\wt\sigma_q^{b_\ell}$ if $(\frac{\alpha_\ell}\ell)=(-1)^{\frac{b_\ell}2}$, and
$\wt{\text{Fr}}_\ell=\wt\sigma_q^{b_\ell}\varepsilon$ if $(\frac{\alpha_\ell}\ell)=(-1)^{\frac{b_\ell}2+1}$.
So when $\ell\equiv 1\mod 4$ and $\ell\equiv g^{b_\ell}\mod q$ with $2\mid b_\ell$, we have
\[
L_\ell(\wt\rho_j,s)=\begin{cases}
(1-\zeta_{2(q-1)}^{jb_\ell}\ell^{-s})^{-2}&\text{ if }(\frac{\alpha_\ell}\ell)=(-1)^{\frac{b_\ell}2}\\
(1+\zeta_{2(q-1)}^{jb_\ell}\ell^{-s})^{-2}&\text{ if }(\frac{\alpha_\ell}\ell)=(-1)^{\frac{b_\ell}2+1}.
\end{cases}
\]

Next we compute the local factors $L_2(\wt\rho_j,s)$ and $L_q(\wt\rho_j,s)$. When $(\frac 2q)=-1$,
we know from last section that 2 is ramified in $\wt K/K$. So $L_2(\wt\rho_j,s)=1$ in this case. Now
we assume $(\frac 2q)=1$. Since $I_2=<\sigma_{-1}>$ and 2 is unramifed in $\wt K/K$, we have
$\wt I_2=<\wt\sigma_{-1}>$ or $\wt I_2=<\wt\sigma_{-1}\varepsilon>$.
The matrixes $I+\wt\rho_j(\wt\sigma_{-1})$ and $I+\wt\rho_j(\wt\sigma_{-1}\varepsilon)$ have rank 1.
So $V^{\wt I_2}$ has dimension 1. Write Fr$_2=\sigma_2^{b_2}$ with $2\mid b_2$. As in last section,
we always assume $b_2\equiv 2\mod 4$ if $q\equiv 3\mod 4$. Recall that $P_0$ be the set of the prime
numbers of the form $A^2+64B^2$ with $A,B\in\Bbb Z$. Since
$\wt\rho_j(\wt{\sigma}_2^{b_2})=\zeta_{2(q-1)}^{jb_2}I$, by Prop.4.3 we have
\[
L_2(\wt\rho_j,s)=\begin{cases}1-\zeta_{2(q-1)}^{jb_2}2^{-s}&\mathrm{if}~~q\not\in P_0,16\nmid q^*-1,
~\mathrm{or}~q\in P_0,16\mid q^*-1\\
1+\zeta_{2(q-1)}^{jb_2}2^{-s} &\mathrm{if}~~q\in P_0,16\nmid q^*-1,
~\mathrm{or}~q\not\in P_0,16\mid q^*-1.
\end{cases}
\]

When $q\equiv 3\mod 4$, we know that $q$ is ramified in $\wt K/K$. So $L_q(\wt\rho_j,s)=1$ for
odd $j$ in this case. Assume $q\equiv 1\mod 4$. Since $I_q=<\sigma_q>$ and $q$ is unramifed in
$\wt K/K$, we have $\wt I_q=<\wt\sigma_q>$ or $\wt I_2=<\wt\sigma_q\varepsilon>$. Thus
$V^{\wt I_q}=0$  if $j\ne 0$, and $V^{\wt I_q}$ has dimension 1 if $j=0$.

The Frobenious map Fr$_q$ of $q$ in $G$ modulo $I_q$ is the identity map. So $\wt{\mathrm{Fr}}_q=1$ or
$\varepsilon$. In [Sect.5, 7] we have showed that $q$ splits in $\wt K/K$ if $q\equiv 1\mod 8$ and is inertia
if $q\equiv 5\mod 8$. So $\wt{\mathrm{Fr}}_2=1$ if $q\equiv 1\mod 8$ and $\wt{\mathrm{Fr}}_2=\varepsilon$
if $q\equiv 5\mod 8$. Thus we get
\[
L_q(\wt\rho_j,s)=\begin{cases}
1&\text{ if }j\ne 0\\
1-q^{-s}&\text{ if }j=0,~q\equiv 1\mod 8\\
1+q^{-s}&\text{ if }j=0,~q\equiv 5\mod 8.
\end{cases}\]

We have computed all the local factors. So we have
\begin{equation}\label{eq2p3}\begin{split}
L(\wt\rho_j,s)=&(1-u_q\zeta_{2(q-1)}^{jb_2}2^{-s})^{-1}(1-(-1)^{\frac{q-1}4}q^{-s})^{-n_j}\\
&\times\prod_{\ell\equiv 1,~2\nmid b_\ell\text{ or } \ell\equiv 3,~2\mid b_\ell}
(1-\zeta_{q-1}^{jb_\ell}\ell^{-2s})^{-1}\\
&\times\prod_{\ell\equiv 3,~2\nmid b_\ell}(1+\zeta_{q-1}^{jb_\ell}\ell^{-2s})^{-1}
\prod_{\ell\equiv 1,~2\mid b_\ell}(1-u_\ell\zeta_{2(q-1)}^{jb_\ell}\ell^{-s})^{-2},
\end{split}
\end{equation}
where $u_q=1$ if $q\not\in P_0,16\nmid q^*-1,
~\mathrm{or}~q\in P_0,16\mid q^*-1$ and $u_q=-1$ otherwise; $n_j=0$ if $j\ne 0$ and $n_0=1$;
and $u_\ell=(\frac{\alpha_\ell}\ell)(-1)^{\frac{b_\ell}2}$. Here in the products,
$"\equiv"$ means the congruence modulo 4.

\begin{thm} Except for the Dirichlet $L$-functions, all Artin $L$-functions of the
Galois extension $\wt K/\Bbb Q$ are explicitly given by (5.1) in the case $q=2$ and
by (5.2) in the case $q$ is odd, where in (5.2) $0\le j<q-1,~2\mid j$ if $q\equiv 1\mod 4$
and $0\le j<q-1,~2\nmid j$ if $q\equiv 3\mod 4$.
\end{thm}

\noindent {\bf 5.3 A formula.} Let $\zeta_{\wt K}(s)$ and $\zeta_K(s)$ be the Dedekind zeta functions of
$\wt K$ and $K$ respectively. By Artin's formula of the decomposition of Dedekind zeta functions,
we have
\[
\frac{\zeta_{\wt K}(s)}{\zeta_K(s)}=\prod_{\wt\rho_j}\prod_{\ell:\text{ prime}}
L_\ell(\wt\rho_j,s)^2,
\]
where $\wt\rho_j$ runs over all 2-dimensional irreducible representations of $\wt G$.
When $q=2$, there is only one 2-dimensional irreducible representation of $\wt G$. So the
square of (5.1) gives the formula. When $q$ is odd, by computing
$\prod_{\wt\rho_j}L_\ell(\wt\rho_j,s)$, we get the following
\begin{cor} For a prime number $\ell\ne q$, let $f_\ell=\frac{q-1}{\gcd(b_\ell,q-1)}$ be the order
of $\ell\mod q$ and let $g_\ell=\gcd(b_\ell,q-1)=\frac{q-1}{f_\ell}$. If $q\equiv 1\mod 4$, we have
\[\begin{split}
\frac{\zeta_{\wt K}(s)}{\zeta_K(s)}=&(1-u_q^{f_2}2^{-f_2s})^{-g_2}(1-(-1)^{\frac{q-1}4}q^{-s})^{-2}
\prod_{\ell\equiv 1,~2\nmid b_\ell\text{ or } \ell\equiv 3}(1-\ell^{-f_\ell s})^{-2g_\ell}\\
&\times\prod_{\ell\equiv 1,~2\mid b_\ell}(1-u_\ell^{f_\ell}\ell^{-f_\ell s})^{-2g_\ell},
\end{split}\]
and if $q\equiv 3\mod 4$, we have
\[\begin{split}
\frac{\zeta_{\wt K}(s)}{\zeta_K(s)}=&(1+u_q^{f_2}2^{-f_2s})^{-g_2}
\prod_{\ell\equiv 1,~2\nmid b_\ell}(1+\ell^{-f_\ell s})^{-2g_\ell}
\prod_{\ell\equiv 3}(1-\ell^{-2f_\ell s})^{-g_\ell}\\
&\times\prod_{\ell\equiv 1,~2\mid b_\ell}(1+u_\ell^{f_\ell}\ell^{-f_\ell s})^{-2g_\ell},
\end{split}\]
where $u_q$ and $u_\ell$ are as above.
\end{cor}


\begin{thebibliography}{9}

\bibitem {An}
G. W. Anderson, \emph{Kronecker-Weber plus epsilon}, { Duke Math. J.}
{\bf 114}, 439-475 (2002).

\bibitem {DP}
M. Daberkow, M. Pohst, \emph{On integral bases in relative quadratic extensions},
{ Math. of Computation} {\bf 65}, 319-329(1996).

\bibitem{Das}
P. Das, \emph{Algebraic Gamma monomials and double coverings of
cyclotomic fields}, Trans. Amer. Math. Soc. \textbf{352}, 3557--3594
(2000).

\bibitem {HY}
Y. Hu, L. Yin, \emph{Arithmetic and representations of a non-abelian Galois number field},
{Preprint}, 2006.

\bibitem {IR}
K. Ireland, M. Rosen, \emph{A Classical Introduction to Modern Number Theory}, GTM 84,
{Springer-Verlag}, 1982.

\bibitem{Se}
J-P. Serre, \emph{Linear Representations of finite groups}, GTM 42, Springer-Verlag,
New York Inc. 1977.

\bibitem {YCZ}
L. Yin and C. Zhang, \emph{Arithmetic of quasi-cyclotomic fields}, {Preprint}, 2006.

\bibitem {YQZ}
L. Yin and Q. Zhang, \emph{All double coverings of cyclotomic
fields}, {Math. Zeit.} {\bf 253}, 479-488(2006).

\end{thebibliography}
\end{document}